\documentclass[smallextended]{svjour3}
\usepackage{amsmath,graphicx,indentfirst,enumitem,amssymb,
tikz,color,geometry,url}
\usepackage{tikz-network}
\usepackage[subsection]{algorithm}
\usepackage{color}
\usepackage{amsfonts}
\usepackage{algorithm}
\usepackage{algpseudocode}
\usetikzlibrary{positioning,calc,arrows,automata}
\newcommand{\R}{\mathbb{R}}

  \newcommand{\one}{\mathbf{1}}
\date{}
\begin{document}
\author{Silvia Noschese \and Lothar Reichel}
\institute{Silvia Noschese\at
Dipartimento di Matematica\\ 
SAPIENZA Universit\`a di Roma\\
P.le Aldo Moro 5, 00185 Roma, Italy\\
\email{noschese@mat.uniroma1.it} {\it (corresponding author)}\\
Lothar Reichel\at 
Department of Mathematical Sciences\\ 
Kent State University\\
Kent, OH 44242, USA\\
\email{reichel@math.kent.edu}
}
\title{Enhancing multiplex global efficiency}
\maketitle

\begin{abstract}
Modeling complex systems that consist of different types of objects leads to multilayer 
networks, in which vertices are connected by both inter-layer and intra-layer edges. In 
this paper, we investigate multiplex networks, in which vertices in different layers are 
identified with each other, and the only inter-layer edges are those that connect a vertex
with its copy in other layers. Let the third-order adjacency tensor 
$\mathcal{A}\in\R^{N\times N\times L}$ and the parameter $\gamma\geq 0$, which is 
associated with the ease of communication between layers, represent a multiplex network with
$N$ vertices and $L$ layers. To measure the ease of communication in a multiplex network, 
we focus on the average inverse geodesic length, which we refer to as the multiplex global 
efficiency $e_\mathcal{A}(\gamma)$ by means of the multiplex path length matrix
$P\in\R^{N\times N}$. This paper generalizes the approach proposed in \cite{NR23} for 
single-layer networks. We describe an algorithm based on min-plus matrix multiplication 
to construct $P$, as well as variants $P^K$ that only take into account multiplex paths 
made up of at most $K$ intra-layer edges. These matrices are applied to detect redundant
edges and to determine non-decreasing lower bounds $e_\mathcal{A}^K(\gamma)$ for
$e_\mathcal{A}(\gamma)$, for $K=1,2,\dots,N-2$. 
 Finally, the sensitivity of $e_\mathcal{A}^K(\gamma)$ to changes 
of the entries of the adjacency tensor $\mathcal{A}$ is investigated to determine edges 
that should be strengthened to enhance the multiplex global efficiency the most.
\end{abstract}



\section{Introduction}\label{s1}
Multilayer networks consist of different kinds of edges and possibly different types of 
vertices. This kind of networks arise when one seeks to model a complex system that 
contains connections and objects with different properties; see, e.g., \cite{M1,M2} for an
overview on multilayer networks. In the particular case of multiplex networks, vertices in
different layers are identified with each other, i.e., every vertex in some layer has a 
copy in all other layers and is connected to them. The 
only inter-layer edges are those that connect instances of the same vertex in different layers. For 
instance, when modeling an urban public transportation network made up of metro and bus 
connections, the metro routes and bus routes define intra-layer edges in different layers
and the locations of the metro stations and bus stops define vertices with diverse 
properties; the cost associated with each intra-layer edge accounts for the time needed to
travel from one location to another, whereas the cost $\gamma\geq 0$, which is determined 
by the average amount of time spent, is associated with each transfer between a metro
station and an adjacent bus stop. This gives rise to an inter-layer (undirected) edge 
along which travelers walk. In the context of models for urban public transportation, we 
remark that a more general type of multiplex network where vertices may not be connected 
between all layers has been recently considered in \cite{BS2}. However, in this paper we 
will deal exclusively with the type of multiplex described above.
 
How efficiently communication between the vertices flows through a multiplex can be 
measured with the aid of the path length matrix associated with the network. Single layer 
shortest paths are made up of edges within one layer, whereas multiplex shortest paths may 
make use of inter-layer edges to move between layers. Note that in a multiplex, in which 
intra-layer edge weights are proportional to the distance between the vertices that
the edge connects, or are proportional to the cost of traveling along an edge, the length 
of a path should take into account both the cost of traversing intra-layer edges, i.e., 
the sum of the relevant weights, and the number of movements between layers 
multiplied by $\gamma$. 

Let us introduce some notation and definitions that will be used throughout this paper. A 
multiplex network may be represented by $L$ graphs that share the same set of vertices 
$V_N=\{v_1,v_2,\dots,v_N\}$. The (possibly weighted and/or directed) graph for layer 
$\ell$ is associated with a non-negative intra-layer adjacency matrix 
$A^{(\ell)}=[a_{ij}^{(\ell)}]_{i,j=1,2,\dots,N}\in\R^{N\times N}$, where
$\ell=1,2,\dots,L$. Alternatively, a multiplex network may be represented by a 
non-negative third-order adjacency tensor $\mathcal{A}=[a_{ij}^{(\ell)}]_{i,j=1,2,\dots,N,
\;\ell=1,2,\dots,L }\in\R^{N\times N\times L}$, where $a_{ij}^{(\ell)}>0$ is the weight of
the edge pointing from vertex $v_i$ to vertex $v_j$ in layer $\ell$ (if such an edge 
exists), and $a_{ij}^{(\ell)}=0$ if there is no edge from $v_i$ to $v_j$ in layer $\ell$. 
The graph is assumed to be simple, i.e., it has at most one edge between any two vertices 
and no edge starts and ends at the same vertex. We remark that De Domenico et al. 
\cite{DSOGA} introduced the supra-adjacency matrix $B\in\R^{NL\times NL}$ associated with 
the multiplex, which has the diagonal blocks $A^{(\ell)}$, $\ell=1,2,\dots,L$, and every 
$N\times N$ off-diagonal block is a multiple of the identity matrix, i.e., 
$\gamma I_N\in\R^{N\times N}$ if $\gamma>0$. As mentioned above, the parameter 
$\gamma\geq 0$ represents the average cost of moving from one layer to another. This 
yields the matrix 
\begin{equation}\label{suprad}
B:=B(\gamma)={\rm blkdiag}[A^{(1)},A^{(2)},\dots,A^{(L)}]+
\gamma(\one_L\one_L^T\otimes I_N-I_{NL}),
\end{equation}
where $\one_n$ denotes the $n$-dimensional vector of all ones and $\otimes$ stands for the
Kronecker product; see \cite{BS,DSOGA}.

To measure the ease of communication between the vertices in a multiplex, we compute the 
average inverse geodesic length of the multiplex. To this end, we need to construct the 
multiplex path length matrix  $P=[p_{ij}]_{i,j=1,2,\dots,N}\in\R^{N\times N}$, whose entry
$p_{ij}$ is the length of the shortest paths from vertex $v_i$ to vertex $v_j$, where the
length is determined by the edge weights; if there is no path between these vertices, then 
$p_{ij}=\infty$. To limit the computational cost of this approach, we are interested in 
determining paths that use at most $K$ edges for some $1\leq K<N$. By means of the 
multiplex $K$-path length matrix  
$P^{K}=[p_{ij}^K]_{i,j=1,2,\dots,N}\in\R^{N\times N}$, one can compute the average inverse 
$K$-geodesic length, $e_\mathcal{A}^K(\gamma)$, and in this way determine a lower bound 
for the {\it multiplex global efficiency},
$e_\mathcal{A}(\gamma):=e_\mathcal{A}^{N-1}(\gamma)$. We note that for a variety of 
multiplex networks $P^K=P^{N-1}$ for some $1\leq K\ll N-1$. This is illustrated by 
computed examples presented in this paper. 

It is often desirable to be able to assess the sensitivity of the multiplex global 
efficiency to changes in the edge weights. For instance, if the vertices represent cities 
and the edges represent roads between the cities, with edge weights proportional to the 
amount of traffic on each road, then one may be interested in which road(s) should be 
widened or made narrower to increase or reduce, respectively, communication in the 
multiplex network the most. Applications of our approach include city planning and 
information transmission. To enhance communication by using information given by $P^K$, 
for some $1\leq K<N$, we investigate the sensitivity of $e_\mathcal{A}^K(\gamma)$ to 
changes of the entries of the adjacency tensor $\mathcal{A}$  by studying suitable vertex 
centrality measures or by applying the Perron-Frobenius theory to the ``reciprocal''
multiplex $K$-path length matrix $P^K_{-1}$, whose off-diagonal entries are the 
reciprocals of $p_{ij}^K$. Hence, the matrix $P^K_{-1}$ is nonnegative, irreducible if the
multiplex is connected, and often sparse if $K\ll N$. This way, we can determine edges that 
should be strengthened in order to increase the multiplex global efficiency the most. A 
related approach for single-layer networks is described in \cite{NR23}.

The situation of redundant edges also can be analyzed by means of the information given by
matrices $P^K$. We say that an intra-layer edge is redundant if it is convenient to follow
an alternative path to get from its first vertex to its last vertex. However, one observes 
that in cases of random attacks on or failures of the network, having redundant edges may 
be useful for protecting the network \cite{cnab021}. Furthermore, if it is equally 
convenient to traverse a given intra-layer edge or follow an alternative path, then such 
an edge may be profitably used in case of  bottlenecks (e.g., in the event of a highway 
affected by an accident or exceptional much traffic).

This paper is organized as follows: In Section \ref{s2} we present an algorithm based on 
min-plus multiplication that constructs the multiplex path length matrix. Section \ref{s3}
is concerned with the issue of determining redundant edges in a multiplex network. In 
Section \ref{s4}, we measure the multiplex global efficiency and its estimates that easily
can be computed by means of the multiplex $K$-path length matrix. Section \ref{s5} 
presents algorithms for determining which edge weight should be changed to boost global 
efficiency the most. Changing an edge weight may entail widening streets, decreasing 
travel times on a highway by increasing the travel speed, or decreasing the waiting time 
for trams on a route by increasing the number of trams. Finally, numerical tests for 
multiplex networks are reported in Section \ref{s6} and concluding remarks can be found in
Section \ref{s7}.

\section{The multiplex path length matrix}\label{s2}
To construct the path length matrix associated with the given multiplex network, we will 
make use of {\it min-plus matrix multiplication}, i.e., matrix multiplication in the 
tropical algebra \cite{L}:
\begin{equation*}
C=A\star B: \qquad c_{ij}=\min_{h=1,2,\ldots,n} \{a_{ih}+b_{hj}\},\qquad 1\leq i,j\leq N,
\end{equation*}
with $A=[a_{ij}]_{i,j=1}^N$, $B=[b_{ij}]_{i,j=1}^N$, and
$C=[c_{ij}]_{i,j=1}^N\in\R^{N\times N}$.

The first step consists of setting the vanishing off-diagonal entries of 
$A^{(\ell)}=[a_{ij}^{(\ell)}]_{i,j=1,2,\ldots,N}$, to $\infty$ for $\ell=1,2,\dots,L$. This
gives the third-order tensor 
\begin{equation}\label{calP}
\mathcal{P}=[p_{ij}^{(\ell)}]_{i,j=1,2,\dots,N,\;\ell=1,2,\dots,L}\in
\R^{N\times N\times L},
\end{equation}
with
\[
p_{ij}^{(\ell)}=\left\{\begin{array}{cc}
a_{ij}^{(\ell)},& \mbox{~~if~~} a_{ij}^{(\ell)}\ne 0,\\
\infty, & \mbox{otherwise},
\end{array}\right.
\]
for all $i\ne j$ and $\ell$. Moreover, $p_{ii}^{(\ell)}=0$ for all $i$ and $\ell$. We are 
in a position to construct the {\it multiplex $1$-path length matrix} 
\[
P^1=[p_{ij}^1]_{i,j=1}^N \mbox{~~with~~} p_{ij}^1=\min_{\ell=1,2,\dots,L} 
p_{ij}^{(\ell)}.
\]
The entry $p_{ij}^1$, with  $i\ne j$, either represents the length of the shortest path 
from vertex $v_i$ to vertex $v_j$ made up of a single (intra-layer) edge, or equals 
infinity if there is no edge in any layer from vertex $v_i$ to vertex $v_j$.

\subsection{The case $\gamma=0$}
Let $\gamma=0$. Then one can use the algorithm \emph{function PATHLENGTH\_MATRIX} for 
single-layer networks described in \cite{NR23} to determine the multiplex path length 
matrix $P=P^{N-1}$ by constructing min-plus powers of $P^1$. In more detail, for $K>1$, the
{\it min-plus power} $P^K$ of $P^1$ is given by
\begin{equation}\label{PK}
P^{K}=[p^{K}_{ij}]_{i,j=1}^N\in\R^{N\times N}: \quad p^{K}_{ij}=\min_{h=1,2,\ldots,N} 
\{p^{K-1}_{ih}+p^{1}_{hj}\}, \,\,\mbox {if}\,\, i\ne j,\;\; \mbox{and} \,\, p^{K}_{ij} =0, 
\;\; \mbox{otherwise}.
\end{equation}
The matrix $P^{K}$ gives vertex distances using multiplex paths of at most $K$ intra-layer
edges. In detail, the entry $p^K_{ij}$, with $i\neq j$, represents the length of the 
shortest path from $v_i$ to $v_j$ made up of at most $K$ intra-layer edges. The diagonal 
entries of $P^{K}$ are zero by definition. One has $p_{ij}^{K}=\infty$ if every path from 
$v_i$ to $v_j$ is made up of more than $K$ intra-layer edges, or if there is no path from 
$v_i$ to $v_j$. The intra-layer edges of a shortest path do not necessarily belong to the 
same layer.

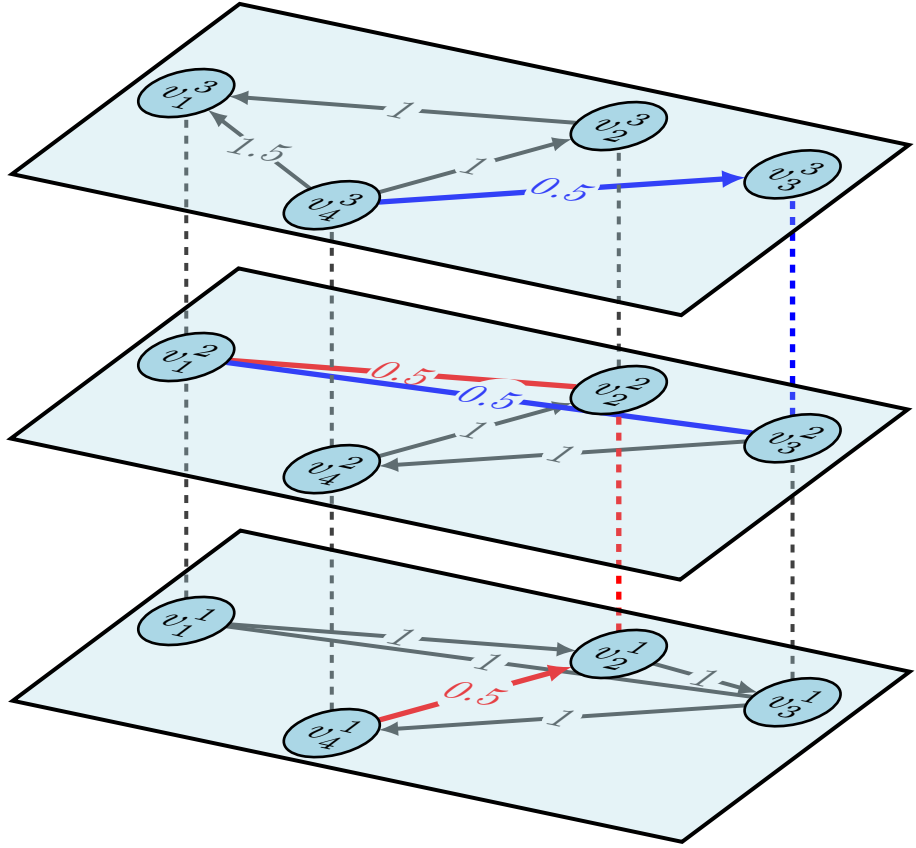
\begin{figure}[ht!]
\centering
\begin{tikzpicture}[multilayer=3d]
\SetLayerDistance{3.5}
\Plane[x=-0.2,y=-3,width=9,height=3.75]
\Plane[x=-3.9,y=1.5,width=9,height=3.75]
\Plane[x=-7.6,y=6.05,width=9,height=3.75]
  \Vertex[x=0.5,y=-1,Math,IdAsLabel,layer=1,fontscale=2,size=1]{v_1^1}
  \Vertex[x=5.5,Math,IdAsLabel,layer=1,fontscale=2,size=1]{v_2^1}
 \Vertex[x=8,y=-0.2,Math,IdAsLabel,layer=1,fontscale=2,size=1]{v_3^1}
 \Vertex[x=3.6,y=-2.4,Math,IdAsLabel,layer=1,fontscale=2,size=1]{v_4^1}
 \Vertex[x=0.5,y=-1,Math,IdAsLabel,layer=2,fontscale=2,size=1]{v_1^2}
  \Vertex[x=5.5,Math,IdAsLabel,layer=2,fontscale=2,size=1]{v_2^2}
 \Vertex[x=8,y=-0.2,Math,IdAsLabel,layer=2,fontscale=2,size=1]{v_3^2}
 \Vertex[x=3.6,y=-2.4,Math,IdAsLabel,layer=2,fontscale=2,size=1]{v_4^2}
 \Vertex[x=0.5,y=-1,Math,IdAsLabel,layer=3,fontscale=2,size=1]{v_1^3}
  \Vertex[x=5.5,Math,IdAsLabel,layer=3,fontscale=2,size=1]{v_2^3}
 \Vertex[x=8,y=-0.2,Math,IdAsLabel,layer=3,fontscale=2,size=1]{v_3^3}
 \Vertex[x=3.6,y=-2.4,Math,IdAsLabel,layer=3,fontscale=2,size=1]{v_4^3}
 
\Edge[label=$1$,Direct,fontscale=2](v_1^1)(v_2^1)
\Edge[label=$1$,Direct,fontscale=2](v_2^1)(v_3^1)
\Edge[label=$1$,fontscale=2](v_1^1)(v_3^1)
\Edge[label=$1$,Direct,fontscale=2](v_3^1)(v_4^1)
\Edge[label=$0.5$,color=red, Direct,fontscale=2,lw=2pt](v_4^1)(v_2^1)

\Edge[label=$0.5$, color=red,fontscale=2,lw=2pt](v_1^2)(v_2^2)
\Edge[label=$1$,Direct,fontscale=2](v_4^2)(v_2^2)
\Edge[label=$0.5$,fontscale=2, color=blue,lw=2pt](v_1^2)(v_3^2)
\Edge[label=$1$,Direct,fontscale=2](v_3^2)(v_4^2)

\Edge[label=$0.5$,Direct,fontscale=2, color=blue,lw=2pt](v_4^3)(v_3^3)
\Edge[label=$1$,Direct,fontscale=2](v_4^3)(v_2^3)
\Edge[label=$1$,Direct,fontscale=2](v_2^3)(v_1^3)
\Edge[label=$1.5$,Direct,fontscale=2](v_4^3)(v_1^3)

\Edge[style=dashed,fontscale=2](v_1^2)(v_1^3)
\Edge[style=dashed,fontscale=2](v_1^1)(v_1^2)  
\Edge[style=dashed,fontscale=2](v_2^2)(v_2^3)  
 \Edge[style=dashed,fontscale=2, color=red,lw=2pt](v_2^1)(v_2^2)  
 \Edge[style=dashed,fontscale=2, color=blue,lw=2pt](v_3^2)(v_3^3) 
\Edge[style=dashed,fontscale=2](v_3^1)(v_3^2)  
\Edge[style=dashed,fontscale=2](v_4^2)(v_4^3) 
\Edge[style=dashed,fontscale=2](v_4^1)(v_4^2)  
\end{tikzpicture}
\caption{Multiplex considered in Examples 1 and 2. When the switching cost $\gamma$ 
satisfies $0\leq\gamma<0.5$, there are two shortest paths from vertex $v_4$ to vertex 
$v_1$: the one shown in red starts at layer $1$ (vertex $v_4^1$) and ends at layer $2$ 
(vertex $v_1^2$) with one switch (through vertices $v_2^1$ and $v_2^2$), whereas the one 
shown in blue starts at layer $3$ (vertex $v_4^3$) and ends at layer $2$ (vertex $v_1^2$) 
with one switch (through vertices $v_3^3$ and  $v_3^2$). 
}\label{fig1}
\end{figure}

\begin{example}\label{ex1}
Three private shuttle services drive between four archaeological sites in one area. In the
corresponding multiplex, the layers represent the companies, the vertices the sites, the 
edges the roads traveled, and the edge weights the average waiting time for a shuttle.
When the shuttles of a company travels along a road at equidistant times by twice as many 
shuttles than the other companies, the average waiting time for a shuttle of this company, 
i.e., the edge weight in the corresponding layer, is $1/2$ instead of $1$ (which is the 
edge weight for the other companies). When a company has a road traversed with $2/3$ of 
the number of shuttles at equidistant times than the other companies, the average waiting 
time, which is the edge weight, for shuttles of this company is $3/2$ instead of $1$. On 
some roads shuttles go back and forth (resulting in undirected edges), on other roads they
only go one way (resulting in directed edges).

We can model the situation described by the supra-adjacency matrix \eqref{suprad},
$$B=\left(
\begin{array}{ccc}
A^{(1)} & \gamma I_4 &\gamma I_4\\  
\gamma I_4 &A^{(2)}& \gamma I_4\\ 
\gamma I_4 &\gamma I_4 &A^{(3)}\\
\end{array}
\right)\in\R^{12\times 12},
$$ 
with the diagonal blocks
$$A^{(1)}=
\left(
\begin{array}{cccc}
0\phantom{00}&1\phantom{00}&1\phantom{00}&0\\
0\phantom{00}&0\phantom{00}&1\phantom{00}&0\\
1\phantom{00}&0\phantom{00}&0\phantom{00}&1\\
0\phantom{00}&0.5&0\phantom{00}&0\\
\end{array}
\right) ,
\;\;\;A^{(2)}=
\left(
\begin{array}{cccc}
0\phantom{00}&0.5&0.5&0\\
0.5&0\phantom{00}&0\phantom{00}&0\\
0.5&0\phantom{00}&0\phantom{00}&1\\
0\phantom{00}&1\phantom{00}&0\phantom{00}&0\\
\end{array}
\right),
\;\;\;A^{(3)}=
\left(
\begin{array}{cccc}
0\phantom{00}&0\phantom{00}&0\phantom{00}&0\\
1\phantom{00}&0\phantom{00}&0\phantom{00}&0\\
0\phantom{00}&0\phantom{00}&0\phantom{00}&0\\
1.5&1\phantom{00}&0.5&0\\
\end{array}
\right);
$$
see Figure \ref{fig1} for a visualization of the associated multiplex.

Assume that the companies pick up tourists from and bring them to the same stops. Hence  
one has $\gamma=0$. The multiplex path length matrix then is given by
$$P=P^2=\left(
\begin{array}{cccc}
0\phantom{00}&0.5&0.5&1.5\\
0.5&0\phantom{00}&1\phantom{00}&2\phantom{00}\\
0.5&1\phantom{00}&0\phantom{00}&1\phantom{00}\\
1\phantom{00}&0.5&0.5&0\phantom{00}\\
\end{array}
\right)\in\R^{4\times 4}.
$$ 
Note that there are two shortest paths from vertex $v_4$ to vertex $v_1$; one is drawn in 
red and the other in blue in Figure \ref{fig1}. Both these paths are made up of two 
intra-layer edges and a (free) layer switch. This means that the traveler, from the site 
associated with $v_4$ to the site associated with $v_1$, may use shuttles from different 
companies.
\end{example}

\subsection{The general case $\gamma\geq0$}
We turn to the situation when $\gamma\geq 0$. When constructing the multiplex $K$-path 
length matrix for $K>1$, the evaluations \eqref{PK} have to be modified because one has to 
include the cost $\gamma$ for each layer switch to the sum of the weights of the 
intra-layer edges of a path. In particular, when determining the length $p_{ij}^{K}$ of a 
shortest path made up of at most $K$ intra-layer edges from vertex $v_i$ to vertex $v_j$, 
for $i,j=1,2,\dots,N$, with $i \ne j$, one has to analyze whether the switching cost 
$\gamma$ is relevant. Specifically, one should consider the layer of the last edge (i.e.,
the intra-layer edge from the penultimate vertex to the last vertex) of any shortest path 
from vertex $v_i$ to vertex $v_h$ made up of at most $K-1$ edges (in case 
$0<p^{K-1}_{ih}<\infty$), and take into account the entries of the third-order tensor 
$\mathcal{P}$ in \eqref{calP} (and not only the entries of $P^1$ as in \eqref{PK}). In 
detail, the off-diagonal entries of the multiplex $K$-path length matrix 
$P^{K}=[p^{K}_{ij}]_{i,j=1}^N\in\R^{N\times N}$ for $K>1$ are computed according to 
\begin{equation}\label{PKgamma}
p^{K}_{ij} =p^{K-1}_{i\bar{h}}+p^{(\bar{\ell})}_{\bar{h}j}+
\gamma \delta_{\bar{h}\bar{\ell}},  \,\, \,\, \,\,\mbox {if}\,\,  i\ne j,
\end{equation}
where 
$$
(\bar{h},\bar{\ell})={\arg\min}_{h=1,2,\ldots,N, \,\ell=1,2,\dots,L} \,
\{p^{K-1}_{ih}+p^{(\ell)}_{hj}+\gamma \delta_{h\ell}\},  
\,\, \,\, \,\,\mbox {if}\,\, i\ne j,
$$
with $\delta_{h\ell}=0$, if one of the following conditions holds:
\begin{itemize}
\item $p^{K-1}_{ih}=0$, i.e., $v_i=v_h$; 
\item there is no path from vertex $v_i$ to vertex $v_h$ made up of at most $K-1$ edges, 
i.e., $p^{K-1}_{ih}=\infty$;
\item $p^{(\ell)}_{hj}=0$, for all $\ell=1,2,\dots,L$, i.e., $v_h=v_j$; 
\item there are no intra-layer edges from vertex $v_h$ to vertex $v_j$, i.e., 
$p^{(\ell)}_{h,j}=\infty$ for all $\ell=1,2,\dots,L$;
\item the intra-layer edge from vertex $v_h$ to vertex $v_j$ with weight 
$p^{({\ell})}_{hj}$ belongs to the same layer $\ell$ of the last edge of a shortest path 
made up of at most $K-1$ edges from vertex $v_i$ to vertex $v_h$ of length $p^{K-1}_{ih}$;
\end{itemize}
and $\delta_{h\ell}=1$ otherwise.  

\begin{example}\label{ex2}
In the model illustrated in Example \ref{ex1}, consider the variant that the shuttle stops
of different companies are located far away; hence one has $\gamma>0$. It easy to see 
that, by means of the evaluations in \eqref{PKgamma}, the multiplex path length matrix is 
given by 
$$P=\left(
\begin{array}{cccc}
0\phantom{00}&0.5&0.5&1.5\\
0.5&0\phantom{00}&1\phantom{00}&2\phantom{00}\\
0.5&1\phantom{00}&0\phantom{00}&1\phantom{00}\\
p_{4,1}(\gamma)&0.5&0.5&0\phantom{00}\\
\end{array}
\right)\in\R^{4\times 4},
$$ 
with $p_{4,1}(\gamma)=\min\{ 1+ \gamma, \,1.5\}$. Indeed, if $\gamma<0.5$, the shortest 
paths from vertex $v_4$ to vertex $v_1$ are the ones drawn in red and in blue in Figure 
\ref{fig1}, whereas if $\gamma = 0.5$ another shortest path from $v_4$ to $v_1$ is given 
by a single intra-layer edge with weight $1.5$ in the third layer. This implies that a 
user may alternatively choose the shuttle of the third company from the site associated 
with $v_4$ to the site associated with $v_1$. Finally, if $\gamma>0.5$, the latter is the
unique shortest path from $v_4$ to $v_1$. One notices that, if $\gamma>0.5$, no shortest 
paths between the sites require layer switches.
\end{example}

\section{Redundant intra-layer edges}\label{s3}
We are interested in which intra-layer edges do not contribute to the network efficiency, 
recalling that an edge is considered redundant if it is convenient to follow an 
alternative route. Note that the triangle inequality holds for the entries of the path 
length matrix $P=[p_{ij}]_{i,j=1}^N\in\R^{N\times N}$. Specifically, 
\begin{equation}\label{triangle}
p_{ij}\leq p_{ih}+p_{hj},\qquad 1\leq i,j\leq N.
\end{equation}
 
\subsection{The case $\gamma=0$}
Given a multiplex with layer switching cost $\gamma=0$, the redundant intra-layer edges 
can be determined by comparing the third-order tensor $\mathcal{P}$ in \eqref{calP} and
the path length matrix $P$. Note that, by definition, one has $0<p_{ij}^{(\ell)}<\infty$ 
if $i\ne j$ and there exists an edge from vertex $v_i$ to vertex $v_j$ in layer $\ell$ 
with weight $p_{ij}^{(\ell)}$. 

In the case that  $p_{ij}^{(\ell)}>p_{ij}>0$, the intra-layer edge from vertex $v_i$ to 
vertex $v_j$ in  layer $\ell$ is redundant. Indeed, the triangle inequality is not 
satisfied by the length of such an intra-layer edge, because there is at least one 
shortest path from vertex $v_i$ to vertex $v_j$, whose length satisfies the triangle 
inequality \eqref{triangle}. Moreover, the redundancy of some intra-layer edges may be 
inferred in advance by comparing $\mathcal{P}$ and a suitable $K$-path length matrix 
$P^K=[p_{ij}^K]_{i,j=1}^N$ with $1\leq K<N-1$. Indeed, if $p_{ij}^{(\ell)}>p_{ij}^K$, then
the intra-layer edge from vertex $v_i$ to vertex $v_j$ in layer $\ell$ is surely 
redundant, because one has $p_{ij}^{K}\geq p_{ij}$ for all $1\leq K\leq N-1$.

\subsection{The general case $\gamma\geq0$}
If the switching cost $\gamma$ in the multiplex is nonnegative, then the redundant 
intra-layer edges can still be determined by comparing the third-order tensor 
$\mathcal{P}$ and the path length matrix $P$ but further analysis is needed. 
Note that the shortest paths from vertex $v_i$ to vertex $v_j$, 
which have length $p_{ij}$, start and arrive at possibly different layers. For $i\ne j$, 
we denote the set containing each layer of the first intra-layer edges (i.e., edges from 
vertex $v_i$) of such shortest paths by $L_{ij}^{(\mathrm{s})}$, and we denote the set
containing each layer of their last intra-layer edges (i.e., to vertex $v_j$) by
$L_{ij}^{(\mathrm{a})}$.

In the case when the edge from vertex $v_i$ to vertex $v_j$ in layer $\ell$ has weight
$p_{ij}^{(\ell)}=p_{ij}$, such an intra-layer edge is a shortest path from vertex $v_i$ to
vertex $v_j$, hence it is surely nonredundant. However, even if $p_{ij}^{(\ell)}$ 
satisfies the inequalities $p_{ij}<p_{ij}^{(\ell)}< p_{ij}+2\gamma$, and $\ell$ does not 
belong to the sets $L_{ij}^{(\mathrm{s})}$ and $L_{ij}^{(\mathrm{a})}$, then the edge from
vertex $v_i$ to vertex $v_j$  in layer $\ell$ is useful, e.g., for a traveler who is at 
the location represented by vertex $v_i$ in layer $\ell$ and has to go to the location 
represented by vertex $v_j$ in layer $\ell$. In fact, the cost of traversing such an
intra-layer edge is less than that of first making a layer switch, then walking a shortest 
path (of length $p_{ij}$), and finally making a second layer switch. 
We therefore say that the edge from vertex $v_i$ to vertex $v_j$ in layer $\ell$ is 
redundant if
\begin{equation}\label{red}
p_{ij}^{(\ell)}>p_{ij}+\gamma\left(\delta_{\ell_{ij}^{(\mathrm{s})}}+
\delta_{\ell_{ij}^{(\mathrm{a})}}\right),
\end{equation}   
with 
\begin{itemize}
\item $\delta_{\ell_{ij}^{(\mathrm{s})}}=1$ if $\ell \notin L_{ij}^{(\mathrm{s})}$, and 
$\delta_{\ell_{ij}^{(\mathrm{s})}}=0$ otherwise, 
\item $\delta_{\ell_{ij}^{(\mathrm{a})}}=1$ if $\ell \notin L_{ij}^{(\mathrm{a})}$, and 
$\delta_{\ell_{ij}^{(\mathrm{a})}}=0$ otherwise,
\end{itemize}
because surely there exists a shorter route in the multiplex, possibly made up of both 
intra-layer edges and  inter-layer edges, the latter ones of cost $\gamma$, from vertex 
$v_i$ to vertex $v_j$ in layer $\ell$. We refer to 
any intra-layer edge that does not satisfy \eqref{red} as a nonredundant edge.
 Moreover, an indication of the redundancy of some 
intra-layer edges may be found in advance by comparing $\mathcal{P}$ and a $K$-path length
matrix $P^K=[p_{ij}^K]_{i,j=1}^N$, with $1\leq K<N-1$. Indeed, for a given intra-layer 
edge from vertex $v_i$ to vertex $v_j$  in layer $\ell$, one has
\[
0<p_{ij}= p_{ij}^{N-1}\leq \dots \leq 
p_{ij}^1\leq p_{ij}^{(\ell) }<\infty.
\]
Thus, such an intra-layer edge is surely redundant if, for a certain $K$, one has
\begin{equation}\label{redK}
p_{ij}^{(\ell)}>p_{ij}^K+\gamma\left(\delta_{\ell_{ij}^{K(\mathrm{s})}}+
\delta_{\ell_{ij}^{K(\mathrm{a})}}\right),
\end{equation} 
with
\begin{itemize}
\item $\delta_{\ell_{ij}^{K(\mathrm{s})}}=1$ if $\ell \notin L_{ij}^{K(\mathrm{s})}$, and 
$\delta_{\ell_{ij}^{K(\mathrm{s})}}=0$ otherwise, 
\item $\delta_{\ell_{ij}^{K(\mathrm{a})}}=1$ if $\ell \notin L_{ij}^{K(\mathrm{a})}$, and 
$\delta_{\ell_{ij}^{K(\mathrm{a})}}=0$ otherwise,
\end{itemize}
where, for $i\ne j$, $L_{ij}^{K(\mathrm{s})}$ denotes the set containing each layer of the 
first intra-layer edges (from vertex $v_i$) of shortest paths made up of at most $K$ 
intra-layer edges and where $L_{ij}^{K(\mathrm{a})}$ is the set containing each layer of 
the last intra-layer edges (to vertex $v_j$). We refer to any intra-layer edge that
does not satisfy \eqref{redK} as a $K$-nonredundant edge.

Note that \eqref{redK} might not be satisfied for all $K<N-1$ by redundant intra-layer 
edges in $\mathcal{P}$: only by constructing $P^{N-1}=P$ it can be excluded that an edge 
is redundant. Indeed, the absence of redundancy is ensured by the fact that the 
inequalities \eqref{red} are not satisfied by any off-diagonal entry of $\cal{P}$.

\begin{example}\label{ex3}
As an illustration of the redundancy of intra-layer edges, we again consider the multiplex
in Example \ref{ex1} with layer switching cost $\gamma\geq 0$. First, let 
$0\leq\gamma<0.5$. Then
$$P^1=\left(
\begin{array}{cccc}
0\phantom{00}&0.5&0.5&\infty\\
0.5&0\phantom{00}&1\phantom{00}&\infty\\
0.5&\infty&0\phantom{00}&1\phantom{00}\\
1.5&0.5&0.5&0\phantom{00}\\
\end{array}
\right),\;\;
P^2=P^3=\left(
\begin{array}{cccc}
0\phantom{00}&0.5&0.5&1.5\\
0.5&0\phantom{00}&1\phantom{00}&2\phantom{00}\\
0.5&1\phantom{00}&0\phantom{00}&1\phantom{00}\\
1+\gamma&0.5&0.5&0\phantom{00}\\
\end{array}
\right)\;\;
$$ 
We notice that 
\[
\begin{array}{lll}
L_{12}^{1(\mathrm{s})}=L_{12}^{1(\mathrm{a})}=\left\{ 2 \right\}, &\,
L_{13}^{1(\mathrm{s})}=L_{13}^{1(\mathrm{a})}=\left\{ 2 \right\}, &\,
L_{14}^{1(\mathrm{s})}=L_{14}^{1(\mathrm{a})}=\emptyset, \\
L_{21}^{1(\mathrm{s})}=L_{21}^{1(\mathrm{a})}=\left\{ 2 \right\}, &\,
L_{23}^{1(\mathrm{s})}=L_{23}^{1(\mathrm{a})}=\left\{ 1 \right\}, &\,
L_{24}^{1(\mathrm{s})}=L_{24}^{1(\mathrm{a})}=\emptyset, \\
L_{31}^{1(\mathrm{s})}=L_{31}^{1(\mathrm{a})}=\left\{ 2 \right\}, & \,
L_{32}^{1(\mathrm{s})}=L_{32}^{1(\mathrm{a})}=\emptyset, &\,
L_{34}^{1(\mathrm{s})}=L_{34}^{1(\mathrm{a})}=\left\{ 1,2 \right\}, \\
L_{41}^{1(\mathrm{s})}=L_{41}^{1(\mathrm{a})}=\left\{ 3 \right\}, &\,
L_{42}^{1(\mathrm{s})}=L_{42}^{1(\mathrm{a})}=\left\{ 1 \right\}, & \,
L_{43}^{1(\mathrm{s})}=L_{43}^{1(\mathrm{a})}=\left\{ 3\right\}.
\end{array}
\]
By checking \eqref{redK} for all $(i,j,\ell)$ and $K=1$, one can see that the edges in the
first layer from $v_1$ to $v_2$, from $v_1$ to $v_3$, and from  $v_3$ to $v_1$ are
redundant, as well as the edge in the second layer from $v_4$ to $v_2$ and the edges in 
the third layer from $v_2$ to $v_1$ and from $v_4$ to $v_2$. Conversely, the redundancy of
the intra-layer edge from $v_4$ to $v_1$ in the third layer is apparent only by looking at
$P^2$ and observing that
$$
L_{41}^{(\mathrm{s})}=L_{41}^{2(\mathrm{s})}=\left\{ 1,3 \right\} \,\, \mathrm{and} \,\
L_{41}^{(\mathrm{a})}=L_{41}^{2(\mathrm{a})}=\left\{ 2 \right\}.
$$
We turn to the situation when $\gamma\geq 0.5$. Then 
$$P^1=\left(
\begin{array}{cccc}
0\phantom{00}&0.5&0.5&\infty\\
0.5&0\phantom{00}&1\phantom{00}&\infty\\
0.5&\infty&0\phantom{00}&1\phantom{00}\\
1.5&0.5&0.5&0\phantom{00}\\
\end{array}
\right),\;\;
P^2=P^3=\left(
\begin{array}{cccc}
0\phantom{00}&0.5&0.5&1.5\\
0.5&0\phantom{00}&1\phantom{00}&2\phantom{00}\\
0.5&1\phantom{00}&0\phantom{00}&1\phantom{00}\\
1.5&0.5&0.5&0\phantom{00}\\
\end{array}
\right),\;\;
$$ 
so that, by \eqref{redK}, one has that the intra-layer edge from $v_4$ to $v_1$ 
in the third layer is nonredundant and that no redundant edges are revealed by only 
looking at $P^2$. 
\end{example}

\section{Multiplex global efficiency}\label{s4}
Analogously to the single-layer case, the \emph{diameter} of a multiplex, represented by a 
third-order adjacency tensor $\mathcal{A}$ and a coefficient $\gamma\geq 0$, can be 
defined as the maximal length $d_{\mathcal{A}}(\gamma)$ of the shortest path between any 
distinct vertices of the multiplex. The diameter provides a measure of how difficult it 
is for the vertices of the network to communicate. Moreover, similarly as in the 
single-layer case \cite{BBV}, the \emph{efficiency} of a path between any two vertices of
a multiplex can be defined as the inverse of the length of the path. As in 
\cite{BBV,NR23}, we refer to the sum $h_i^{\rm out}$ of the efficiencies of all shortest 
paths starting from $v_i$, i.e.,
\[
h_i^{\rm out}=\sum_{j\ne i}p_{ij}^{-1},
\]
as the \emph{harmonic out-centrality} of $v_i$, and the sum $h_j^{\rm in}$ of the 
efficiencies of all shortest paths ending at $v_j$, i.e.,
\[
h_j^{\rm in}=\sum_{i\ne j}p_{ij}^{-1},
\]
as the \emph{harmonic in-centrality} of $v_j$. These measures give a large centrality to 
vertices that have small shortest path distances to/from other vertices of the multiplex. 

If the multiplex is {\it connected}, then the average shortest path efficiency over all 
possible pairs is referred to as the \emph{global efficiency} of the network \cite{BBV}:
\begin{equation}\label{globeff}
e_{\mathcal{A}}(\gamma)=\frac{1}{N(N-1)}{\sum_{i,j\ne i} p_{ij}^{-1}}.
\end{equation}
Note that the measure $e_{\mathcal{A}}(\gamma)$ also is useful when the multiplex has more
than one connected component, because infinite distances do not contribute to the sum 
\eqref{globeff}. Networks with large global efficiency are easy to navigate, a desirable 
property of transportation networks.

As in the single-layer case \cite{NR23}, we introduce the {\it reciprocal $K$-path length 
matrix} 
\[
P^{K}_{-1}=[p^{(K,-1)}_{ij}]_{i,j=1}^N,
\]
which is obtained by replacing the off-diagonal entries of the $K$-path length matrix 
$P^K$, for $1\leq K\leq N-1$, by their reciprocals, i.e., 
\[
p^{(K,-1)}_{ij}=1/p^K_{ij},\qquad 1\leq i,j\leq N,\;\; i\ne j,
\]
where $1/\infty$ is identified with $0$. We also introduce the \emph{harmonic 
$K^{\rm out}$-centrality} of the vertex $v_i$, 
\[
h_{K,i}^{\rm out}=\sum_{j\ne i}p^{(K,-1)}_{ij},
\]
the \emph{harmonic $K^{\rm in}$-centrality} of the vertex $v_j$, 
\[
h_{K,j}^{{\rm in}}=\sum_{i\ne j}p^{(K,-1)}_{ij},
\]
as well as the {\it global $K$-efficiency} of the multiplex
\begin{equation*}
e_{\mathcal{A}}^{K}(\gamma)=\frac{1}{N(N-1)}\sum_{i,j\ne i} p_{ij}^{(K,-1)}=
\frac{1}{N(N-1)}\mathbf{1}_N^TP^{K}_{-1}\mathbf{1}_N, \;\;  \;\; 1\leq K\leq N-1.
\end{equation*}

\begin{example}\label{ex4}
Consider again the multiplex of Example \ref{ex1} represented by $\mathcal{A}$. One has
$d_{\mathcal{A}}(\gamma)=2$ for all $\gamma\geq0$. Table \ref{tab1} reports the global
$K$-efficiency of $\mathcal{A}$ for $K=1$ and $K=2$, and for several values of $\gamma$. 
The table shows the global $K$-efficiency to be independent of $\gamma$ for $K=1$, but 
$e_{\mathcal{A}}(\gamma)=e_{\mathcal{A}}^{2}(\gamma)$ achieves its maximum for $\gamma=0$
and attains its minimum value for all $\gamma\geq 0.5$. Table \ref{tab2} displays the 
harmonic $K^{\rm{in}}$-centrality and the harmonic $K^{\rm out}$-centrality of all 
vertices for $\gamma=0.5$. The table shows that the harmonic 
$K^{\rm in}$-centrality of $v_2$ and $v_4$ and the harmonic $K^{\rm out}$-centrality of 
$v_1$, $v_2$, and $v_3$ increase with $K$. The results mentioned can be expected by 
looking at the matrices $P^1$ and $P^2$. 
We remark that both harmonic $K^{\rm in}$-centrality and harmonic $K^{\rm out}$-centrality 
are independent of $\gamma$ for $K=1$, because the matrix $P^1$ is independent of 
$\gamma$. 
\end{example}

\begin{table}
\centering
\begin{tabular}{l*{2}{|c}}
$\gamma$ & $e_{\mathcal{A}}^{1}(\gamma)$ & $e_{\mathcal{A}}^{2}(\gamma)$ \\
\hline
$0$            & 1.2222  & 1.4306 \\
$0.25$            & 1.2222 & 1.4139 \\
$0.50$            & 1.2222& 1.4028 \\
$0.75$            & 1.2222 & 1.4028\\
$1$            & 1.2222& 1.4028 \\
\end{tabular}
\caption{Example \ref{ex4}. Global $K$-efficiency for $K=1$ and $K=2$ of the multiplex 
depicted in Figure \ref{fig1} for $\gamma=0:0.25:1$.}\label{tab1}
\end{table}

\begin{table}[ht]
\centering
\begin{tabular}{l*{2}{|c}}
$i$ & $h_{1,i}^{{\rm in}}(0.5)$ &  $h_{2,i}^{{\rm in}}(0.5)$ \\
\hline
$1$           & 4.6667 & 4.6667   \\
$2$          & 4.0000  & 5.0000\\
$3$          & 5.0000  &  5.0000\\
$4$         & 1.0000  &  2.1667\\
\end{tabular}
\quad\quad\quad
\begin{tabular}{l*{2}{|c}}
$j$           &  $h_{1,j}^{{\rm out}}(0.5)$ &  $h_{2,j}^{{\rm out}}(0.5)$  \\
\hline
$1$           & 4.0000  & 4.6667   \\
$2$          & 3.0000  & 3.5000\\
$3$          & 3.0000  &  4.0000\\
$4$         & 4.6667  &  4.6667\\
\end{tabular}
\caption{Example \ref{ex4}. Harmonic $K^{\rm in}$-centrality (left table) and harmonic 
$K^{\rm out}$-centrality (right table) for $K=1$ and $K=2$ for all vertices of the 
multiplex depicted in Figure \ref{fig1} for $\gamma= 0.5$.}\label{tab2}
\end{table}

\section{Enhancing global  efficiency} \label{s5}
Section \ref{s3} addressed the issue of determining intra-layer edges that can be 
removed without affecting the global efficiency of the multiplex. This section considers
the converse problem of determining which intra-layer edges contribute the most to the 
network efficiency. 

Let  $\mathcal{S}_+\subset \R^{N\times N}$ be the subspace formed by 
the matrices with the same zero-structure of the aggregated adjacency matrix 
\[
A_+=\sum_{\ell=1}^L A^{(\ell)}.
\]
Let $M|\mathcal{S}_+$ denote the ``projection'' of the matrix $M$ onto $\mathcal{S}_+$, 
i.e., $M|\mathcal{S}_+$ is obtained by setting all the entries of $M$ to $0$ that are $0$ 
in $A_+$. To increase the global efficiency of a single-layer network as much as possible,
by considering paths that connect two vertices, say $v_h$ and $v_k$, it has been shown in
\cite{NR23} that one can apply two different approaches to decide which connection(s) to 
strengthen. Consider only paths made up of $K$ edges. Refining the argument in \cite{NR23}
and adapting it to the multiplex case, one can choose to exploit the information given by
\begin{enumerate}[label=(\roman*)]
\item  $K^{\rm in}$- and $K^{\rm out}$-centralities, $\mathbf{h}_K^{\rm in}$ and 
$\mathbf{h}_K^{\rm out}$, by strengthening existing connections corresponding to the 
largest entry of $(\mathbf{h}_K^{\rm in}\mathbf{h}_K^{\rm out\,T})|\mathcal{S}_+$;
\item the left and right Perron vectors, $\mathbf{y}_K$ and $\mathbf{x}_K$, associated 
with the reciprocal $K$-path length matrix, by strengthening existing connections 
corresponding to the largest entry of $(\mathbf{y}_K\mathbf{x}_K^T)|\mathcal{S}_+$.
\end{enumerate}

\begin{remark}\label{r1}
Our purpose here is to investigate how one can enhance the global efficiency of the
multiplex represented by a given adjacency tensor $\mathcal{A}\in \R^{N\times N\times L}$
by considering only paths made up of $K$ edges, i.e., by only using information given by 
$P^K\in \R^{N\times N}$ for some $1\leq K<N$.  We have observed that the global efficiency 
is insensitive to changes of any entry $a_{ij}^{(\ell)}$ of $\mathcal{A}$ such that 
$0<a_{ij}^{(\ell)}=p_{ij}^{(\ell)}$ represents the weight of a redundant intra-layer edge. 
Notice that by means of the information given by $P^K$, one knows that if  
$p_{ij}^{(\ell)}$ satisfies the inequality in \eqref{redK}, then such intra-layer edge is 
surely redundant. We recall that it can be excluded that a $K$-nonredundant edge be 
redundant only when $P^K=P$.
\end{remark}

Recalling that strengthening is achieved by \emph{decreasing} appropriate weights, we 
``perturb'' the supra-adjacency matrix as follows:
\begin{equation}\label{Btilde1}
\widetilde{B}:=\widetilde{B}(\gamma)=
\rm {blkdiag}[\tilde{A}^{(1)},\tilde{A}^{(2)},\dots,\tilde{A}^{(L)}]+
\gamma(\one_L\one_L^T\otimes I_N-I_{NL}),
\end{equation}
where
\begin{equation}\label{Aelltilde1}
\tilde{A}^{(\ell)}= A^{(\ell)}+\alpha_{h,k}^{(\ell)}\mathbf{e}_{h}\mathbf{e}_{k}^T, 
 \;\;\; \mbox{with}\;\;\alpha_{h,k}^{(\ell)}= - a_{h,k}^{(\ell)}/2, 
% \;\;\;\; \bl{\ell=1,2,\dots,L},
\end{equation}
with the index pairs $(h,k)$ determined by one of the above procedures $(i)$ or $(ii)$ and
$\ell$ such that $a_{h,k}^{(\ell)}>0$  is the weight of a $K$-nonredundant edge (cf.
Remark \ref{r1}).
Here, $\mathbf{e}_i \in \R^N$ denotes the  vector with all zero 
entries except for the $i$th entry, which is one.

If the graph associated with $A^{(\ell)}$ is undirected, then $ \tilde{A}^{(\ell)}$ is
defined as
\begin{equation}\label{Aelltilde2}
\tilde{A}^{(\ell)}=A^{(\ell)}+\alpha_{h,k}^{(\ell)}(\mathbf{e}_{h}\mathbf{e}_{k}^T +
\mathbf{e}_{h}\mathbf{e}_{k}^T).
\end{equation}

\subsection{Harmonic centralities versus eigenvector centralities}
The first approach outlined above is easy to explain: One strengthens
any $K$-nonredundant edge from a
vertex that quickly collects information (i.e., a vertex with the highest harmonic 
$K^{\rm in}$-centrality) to a vertex that quickly broadcasts the information (i.e., a 
vertex with the highest harmonic $K^{\rm out}$-centrality). That is one strengthens 
intra-layer edges from $v_h$ to $v_k$ with
\begin{equation}\label{HK}
(h,k):  h^{\rm in}_{K,h}\,h^{\rm out}_{K,k}=
(\mathbf{h}^{\rm in}_K\mathbf{h}^{\rm out\,T}_K)_{h, k}=
\max_{\substack{i,j=1,2,\ldots,N\\A_+(i,j)>0}} 
(\mathbf{h}^{\rm in}_K\mathbf{h}^{\rm out\,T}_K)_{i,j}.
\end{equation}

The second approach is based on Perron-Frobenius theory. Assume the reciprocal path length
matrix $P^K_{-1}$ is irreducible (i.e., the multiplex is connected). Then its left and
right Perron vectors $\mathbf{y}_K=(y_{K,1},y_{K,2},\ldots,y_{K,N})^T$ and 
$\mathbf{x}_K=(x_{K,1},x_{K,2},\ldots,x_{K,N})^T$, respectively, of unit norm and with 
positive entries are unique. Let $\rho_K$ denote the Perron root. The Perron vectors 
determine the Wilkinson perturbation
\[
W_K=\mathbf{y}_K\mathbf{x}_K^T;
\]
see \cite[Section 2]{W}. Using the technique in \cite{NR,sm1}, in order to induce the 
maximal perturbation in $\rho_K$, one chooses the index pair $(h,k)$ such that $W_K(h,k)$ 
is the largest entry of $W_K$ and $A_+(h,k)>0$, i.e., the indices of the largest entry of 
the Wilkinson perturbation projected onto the zero-structure of $A_+$; see, e.g.,
\cite{NP06} for further details. Thus,
\begin{equation}\label{WK}
(h,k):  x_{K,h}\,y_{K,k}=(\mathbf{y}_K\mathbf{x}_K^T)_{h, k}=
\max_{\substack{i,j=1,2,\ldots,N\\A_+(i,j)>0}} (W_K)_{i,j}.
\end{equation}
We expect the global $K$-efficiency to increase the most when decreasing the weights
of the $K$-nonredundant edges that make the Perron root $\rho_K$ change the most.  

Note that the difference in \eqref{HK} and \eqref{WK} is analogous to the difference 
between considering the vertex with the largest degree the most important vertex, and 
considering the vertex with maximal eigenvector centrality the most important vertex. Both
approaches maximize lower bounds for the global $K$-efficiency of the multiplex. In fact, 
the $1$-norm of $\mathbf{h}^{\rm in}_K$, which coincides with the $1$-norm of 
$\mathbf{h}^{\rm out}_K$, is exactly the sum in the numerator of the global 
$K$-efficiency, while the $\infty$-norm of $\mathbf{h}^{\rm in}_K$ and 
$\mathbf{h}^{\rm out}_K$ are in turn the $1$-norm and the $\infty$-norm of $P^{K}_{-1}$,
respectively. Therefore, one has
\[
N(N-1)\,e^K_{\mathcal{A}}=\|\mathbf{h}^{\rm in}_K\|_1= \|\mathbf{h}^{\rm out}_K\|_1\geq 
\max(\|\mathbf{h}^{\rm in}_K\|_{\infty}, \|\mathbf{h}^{\rm out}_K\|_{\infty})=
\max(\|P^{K}_{-1}\|_{1},\|P^{K}_{-1}\|_{\infty}) \geq\rho_K.
\]

While we determine the vertex importance in the present paper, it may also be interesting 
to calculate the edge importance. An approach for single-layer networks that is based on
the use of the line graph for the network is described in \cite{DLCMR}. The computations
are somewhat complicated and an extension to multilayer networks is outside the scope of
the present paper. An approach to approximate the edge importance for a single-layer 
network by using the vertex importance is described by Arrigo and Benzi \cite{AB}. This
approach does not always identify the most important edges correctly, but the computations
are simple. We note that in our approach a small edge weight makes an edge important, 
while in \cite{AB,DLCMR} a large edge weight makes an edge important. 

\begin{example}\label{ex5}
We apply the above procedures to the multiplex of Example \ref{ex1}. For $0\leq\gamma<0.5$ 
both \eqref{HK} and \eqref{WK} yield $(h,k)=(3,4)$ for both $K=1$ and $K=2$. As for the 
diagonal blocks of  the perturbed supra-adjacency matrix in \eqref{Btilde1},
one has
$$ \tilde{A}^{(1)}=
\left(
\begin{array}{cccc}
0\phantom{00}&1\phantom{00}&1\phantom{00}&0\phantom{00}\\
0\phantom{00}&0\phantom{00}&1\phantom{00}&0\phantom{00}\\
1\phantom{00}&0\phantom{00}&0\phantom{00}&0.5\\
0\phantom{00}&0.5&0\phantom{00}&0\phantom{00}\\
\end{array}
\right),
\;\;\; \tilde{A}^{(2)}=
\left(
\begin{array}{cccc}
0\phantom{00}&0.5&0.5&0\phantom{00}\\
0.5&0\phantom{00}&0\phantom{00}&0\phantom{00}\\
0.5&0\phantom{00}&0\phantom{00}&0.5\\
0\phantom{00}&1\phantom{00}&0\phantom{00}&0\phantom{00}\\
\end{array}
\right),
\;\;\; \tilde{A}^{(3)}={A}^{(3)}.
$$
The same pair $(h,k)$ is obtained by both \eqref{HK} and \eqref{WK} also for 
$\gamma\geq 0.5$ when $K=1$. We report in Table \ref{tab3} (left-hand side table) the 
global $K$-efficiency of
the multiplex represented by the perturbed adjacency tensor $\widetilde{\mathcal{A}}$ for
$K = 1$ and $ K = 2$, and for the same values of $\gamma$ considered in Table \ref{tab1}, 
restricting $\gamma$ to be smaller than $0.5$ for $K=2$.

\begin{table}
\centering
\begin{tabular}{l*{2}{|c}}
$\gamma$ & $e_{\widetilde{\mathcal{A}}}^{1}(\gamma)$ & 
$e_{\widetilde{\mathcal{A}}}^{2}(\gamma)$ \\
\hline
$0$            & 1.3056  & 1.5556 \\
$0.25$            & 1.3056 & 1.5389 \\
$0.50$            & 1.3056& -- \\
$0.75$            & 1.3056 & --\\
$1$            & 1.3056& -- \\
\end{tabular}
\quad \quad \quad
\begin{tabular}{l*{2}{|c}}
$\gamma$ & $e_{\widetilde{\mathcal{A}}}^{1}(\gamma)$ & 
$e_{\widetilde{\mathcal{A}}}^{2}(\gamma)$ \\
\hline
$0$            & -- & -- \\
$0.25$            & -- & -- \\
$0.50$            & --&1.6083 \\
$0.75$            & -- & 1.5972\\
$1$            & --& 1.5972\\
\end{tabular}
\caption{Example \ref{ex5}. Global $K$-efficiency, for $K=1$ and $K=2$ and 
$\gamma=0:0.25:1$, of the multiplex depicted in Figure \ref{fig1} after being perturbed 
differently for different values of $\gamma$. In detail, according to the procedures 
\eqref{HK} and \eqref{WK}, the perturbed multiplex considered in the left-hand side table 
has been obtained by strengthening the intra-layer edges from vertex $v_3$ to vertex 
$v_4$, whereas the values of global efficiency in the right-hand side table can been 
obtained by strengthening either the intra-layer edge from vertex $v_2$ to vertex $v_1$ or
the intra-layer edge from vertex $v_3$ to vertex $v_1$ in the second layer of the original 
multiplex.}\label{tab3}
\end{table}

On the contrary, for $K=2$ and  $\gamma\geq 0.5$, by both \eqref{HK} and \eqref{WK} one 
obtains either $(h,k)=(2,1)$, so that the matrices in \eqref{Aelltilde1} are
\[
\tilde{A}^{(1)}={A}^{(1)}, \;\;\; \tilde{A}^{(2)}=
\left(
\begin{array}{cccc}
0\phantom{000}&0.5&0.5&0\\
0.25&0\phantom{00}&0\phantom{00}&0\\
0.5\phantom{0}&0\phantom{00}&0\phantom{00}&1\\
0\phantom{000}&1\phantom{00}&0\phantom{00}&0\\
\end{array}
\right),
\;\;\; \tilde{A}^{(3)}={A}^{(3)},
\]
or $(h,k)=(3,1)$, in which case
\[
\tilde{A}^{(1)}={A}^{(1)}, \;\;\; \tilde{A}^{(2)}=
\left(
\begin{array}{cccc}
0\phantom{000}&0.5&0.5&0\\
0.5\phantom{0}&0\phantom{00}&0\phantom{00}&0\\
0.25&0\phantom{00}&0\phantom{00}&1\\
0\phantom{000}&1\phantom{00}&0\phantom{00}&0\\
\end{array}
\right),
\;\;\; \tilde{A}^{(3)}={A}^{(3)}.
\]
Notice that these intra-layer edges are part of the shortest paths from vertex $v_4$ to 
vertex $v_1$ when $\gamma<0.5$; cf. Figure \ref{fig1}. The values in Table \ref{tab3} 
(right-hand side table) have been computed by taking into account the multiplex, where the
intra-layer edge from vertex $v_2$ to vertex $v_1$ in the second layer [or equivalently 
the intra-layer edge from vertex $v_3$ to vertex $v_1$ in the second layer] has been 
strengthened. Note that both the intra-layer edge from vertex $v_2$ to vertex $v_1$ in
the third layer and the intra-layer edge from vertex $v_3$ to vertex $v_1$ in the first 
layer are redundant (cf. Example \ref{ex3}) so that their strengthening would be useless.
%One has $e_{\widetilde{\mathcal{A}}}^{2}(0.5)=1.6083$ and 
%$e_{\widetilde{\mathcal{A}}}^{2}(0.75)=e_{\widetilde{\mathcal{A}}}^{2}(1)=1.5972$.
\end{example}

\section{Numerical tests} \label{s6} 
The numerical tests reported in this section have been carried out by using MATLAB R2023a 
on a $3.2$ GHz Intel Core i7 6 core iMac. The Perron root, and the  left and right Perron 
vectors for small to moderately sized networks can easily be evaluated by using the MATLAB 
function {\sf eig}. For large-scale multiplexes, these quantities can be computed by the 
MATLAB function {\sf eigs} or by the two-sided Arnoldi algorithm introduced by Ruhe 
\cite{R} and improved by Zwaan and Hochstenbach \cite{ZH}. 

\subsection{European airlines data set}
The European airlines data set consists of $450$ vertices that represent European airports 
and has $L=37$ layers that represent different airlines operating in Europe. Each edge
represents a flight between airports. There are $3588$ edges, which represent available 
routes. Similarly as in \cite{BS,sm2,TPM}, we set $\gamma=1$ to reflect the effort 
required to change airlines for connecting flights, and we only include the $N=417$ 
vertices of the largest connected component of the network. This component can be 
represented by a third-order tensor $\mathcal{A}\in\R^{N\times N \times L}$, where the 
adjacency matrix of the layer corresponding to a given airline contains $1$ if the 
airline offers a flight between the two corresponding airports, and  $0$ otherwise. 
The network can be downloaded from \cite{B_repository}. 

The multiplex is both undirected and unweighted. Since $\gamma=1$, the length of a path is 
given by the total number of intra-layer and inter-layer edges traversed by the path, and 
the diameter is the maximum number of edges traversed by a shortest path. In this network 
one has $d_{\mathcal A}(1)=9$ and the  path length matrix $P=P^{416}$ is equal to $P^{7}$,
because all the maximal shortest paths are made up of seven intra-layer edges and two 
layer switches. Moreover, four pairs of vertices are connected by shortest paths of 
maximal length: ($v_{413}$, $v_{144}$), ($v_{413}$, $v_{202}$), ($v_{413}$, $v_{316}$), 
and ($v_{413}$, $v_{350}$). The multiplex shows that the Le Mans-Arnage Airport 
($v_{413}$) is poorly connected with the Mehamn Airport ($v_{144}$), Valan Airport 
($v_{202}$),  Berlevag Airport ($v_{316}$), and  Batsfjord Airport ($v_{350}$). Reaching 
these airports requires flights operated by three different airlines and six stopovers. We
can observe that there are no redundant edges in the European airlines network. As for 
the global efficiency, one has  $e_{\mathcal{A}}(1)=e_{\mathcal{A}}^7(1)=0.3477$. Both the
choices \eqref{HK} and \eqref{WK}, with $K=7$, return the pair of vertices 
$(v_{40},v_{15})$. The $3^{rd}$, $9^{th}$, $21^{th}$, and $27^{th}$ layers contain edges 
that connect these vertices. If, according to the procedure in \eqref{Aelltilde2}, one 
changes the entries $a_{15,40}^{(\ell)}$ and $a_{40,15}^{(\ell)}$ for each of
the above listed values of $\ell$, one obtains 
\[
e_{\widetilde{\mathcal{A}}}(1)=e_{\widetilde{\mathcal{A}}}^7(1)=0.3486.
\]
This suggests that the number of flights from the Amsterdam Airport Schiphol (vertex
$v_{15}$) to the Barcelona El Prat Airport (vertex $v_{40}$) operated by EasyJet (layer 
$3$), KLM (layer $9$), Vueling (layer $21$), and Transavia Holland (layer $27$) should be
doubled in order to half the wait time between these flights. Doubling the number of
flights corresponds to halving the weight for the corresponding edge. 

Interestingly, the information provided by the reciprocal path length matrix is the same 
as the one given by  $P^{K}_{-1}$ with $K=2$, because the perturbation that increases the
global $2$-efficiency the most is the same that increases the global efficiency the most; 
cf. Table \ref{tab4}. 

\begin{table}[ht]
\centering
\begin{tabular}{c*{3}{|c}}
$K$           & $(h,k)$ & $e_{\mathcal{A}}^K(1)$ & $e_{\widetilde{\mathcal{A}}}^K(1)$  \\
\hline
7          & $(15,40)$ &$3.476599\cdot 10^{-1}$ & $3.486327\cdot 10^{-1}$    \\
6          & $(15,40)$ & $3.476567\cdot 10^{-1}$  & $3.486295\cdot 10^{-1}$    \\
5          & $(15,40)$ & $3.474249\cdot 10^{-1}$  & $3.483962\cdot 10^{-1}$    \\
4          & $(15,40)$ & $3.441131\cdot 10^{-1}$  &  $3.450480\cdot 10^{-1}$ \\
3          & $(15,40)$ & $3.194896 \cdot 10^{-1}$ &  $3.201297\cdot 10^{-1}$ \\
2          & $(15,40)$ & $1.839298\cdot 10^{-1}$  & $1.840478\cdot 10^{-1}$ \\
1          & $(15,12)$ & $3.404584\cdot 10^{-2}$ & $3.405737\cdot 10^{-2} $  \\
\end{tabular}

\caption{European airlines data set. Indices chosen by the procedures and the global 
$K$-efficiency for both the original multiplex and the perturbed multiplex as in eqs.
\eqref{Btilde1},  \eqref{Aelltilde1}, and \eqref{Aelltilde2} for $K=1,2,\ldots,7$.} 
\label{tab4}
\end{table}

\subsection{The Scotland Yard data set}
This data set has been built from the Scotland Yard board game by the authors
of \cite{BS}. The network can be downloaded from \cite{B_repository}. It consists of 
$N=199$ vertices representing public transportation stops in the city of London and has 
$L=4$ layers that represent different modes of transportation: boat, underground, bus, and
taxi. The $3324$ edges are weighted and undirected. Their weights are determined so that 
all edges in the taxi layer have weight one. A taxi ride is defined as a trip by a taxi 
between two adjacent vertices in the taxi layer; a taxi ride along $k$ edges is considered
$k$ taxi rides. The edge weights in the boat, underground, and bus layers are chosen to be 
equal to the minimal number of taxi rides required to travel between the same vertices. 

We let $\gamma=1$. One has $d_{\mathcal A}(1)=20$. The path length matrix $P=P^{198}$ is 
equal to $P^{20}$. There are four pairs of vertices that are connected by shortest paths 
of length $20$. They are ($v_{175}$, $v_{1}$), ($v_{175}$, $v_{8}$), 
($v_{175}$, $v_{18}$), and ($v_{18}$, $v_{106}$). This suggests that staying at the stop 
$v_{175}$ may be a good choice for Mister X, when he has to reveal his location. Moreover,
all connections with $v_{175}$ are in the taxi layer, which is the layer that leaves room 
for more combinations and wreaks havoc among the players playing Scotland Yard detectives. 
As for the global efficiency, one has 
$e_{\mathcal{A}}(1)=e_{\mathcal{A}}^{20}(1)=0.1665$.

 \begin{table}
\centering
\begin{tabular}{c*{7}{|c}}
       & $K=20$ & $K=19$ & $K=18$ & $K=17$ & $K=16$ & $K=15$& $K=14$   \\
\hline
$e_{\mathcal{A}}^K(1)$& $0.1665$&$0.1665$&$0.1665$&$0.1665$&$0.1665$&$0.1665$&
$0.1664$\\
$(h,k)$ by \eqref{HK}& $(126,114) $ & $(126,114) $ & $(126,114) $ & $(126,114) $ & $(126,114) $ & 
$(126,114) $& $(126,114) $\\
$(h,k)$ by \eqref{WK}& $(126,114) $ & $(126,114) $ & $(126,114) $ & $(126,114) $ & $(126,114) $ & 
$(126,114) $& $(126,114) $  \\
\hline
& $K=13$ & $K=12$ & $K=11$ & $K=10$ & $K=9$ & $K=8$& $K=7$   \\
\hline
$e_{\mathcal{A}}^K(1)$& $0.1663$&$0.1660$&$0.1656$&$0.1647$&$0.1633$&$0.1607$&
$0.1556$\\
$(h,k)$ by \eqref{HK}& $(126,114) $ & $(126,114) $ & $(126,114) $ & $(126,114) $ & $(126,114) $ & 
$(140,126) $& $(140,126)$\\
$(h,k)$ by \eqref{WK}& $(126,114) $ & $(126,114) $ & $(126,114) $ & $ (126,114)$ & $(126,114)$ & 
$(126,114)  $& $(126,114) $  \\
\end{tabular}
\caption{Scotland Yard data set. Global $K$-efficiency and indices chosen by the procedures 
\eqref{HK} and \eqref{WK} for $K=7,8,\ldots,20$.}\label{tab5}
\end{table}

Both the choices \eqref{HK} and \eqref{WK}, starting from $K=9$ and $K=7$, respectively, 
determine the pair of indices $(126,114)$; cf. Table \ref{tab5}. The $4^{th}$ layer (the 
taxi layer) contains an edge that connects these vertices. The global efficiency of the 
multiplex perturbed as in \eqref{Aelltilde2}, with $\alpha_{126,114}^{(4)}=-0.5$, is 
$e_{\widetilde{\mathcal{A}}}(1)=0.1678$. 
The information of interest is that the route connecting the taxi stops represented by 
$v_{126}$ and $v_{114}$ constitutes a potential bottleneck. Therefore, players who play 
Scotland Yard detectives should be at one of the two stops, while the player playing 
Mister X, if it is not possible to stay away from these stops, should play one of his 
``double move'' tokens.
 
If $\gamma=0$, the only redundant edge is the one corresponding to 
$a_{67,111}^{(2)}=a_{111,67}^{(2)}$; it has weight $6$, while $p_{67,111}=p_{111,67}=5$. 
Moreover, a shortest path that connects vertex $v_{67}$ with vertex $v_{111}$ in the third 
layer is made up of three intra-layer edges. Hence, this redundancy can already be 
observed in $P^3$. However, since $\gamma=1$ and the shortest path does not directly 
connect the vertices $v_{67}$ and $v_{111}$ in the second layer as the original edge does, 
a player playing Scotland Yard detective at the underground stop $v_{111}$ (i.e., in the 
second layer), who has to go to the underground stop $v_{67}$, has to transfer to the bus 
layer before following the shortest path and transferring back to the underground layer 
afterwords. Thus, the total length of the route will be larger than the weight of the 
intra-layer edge represented by $a_{111,67}^{(2)}$.

%level 3: 111-100 weight 1 // 100 - 82 weight 2 // 82 - 67 weight 2 
%level 3: 111-100 weight 1 // 100 - 82 weight 2 // 82 - 67 weight 2 

\section{Concluding remarks} \label{s7}
The path length matrix associated with a multiplex represented by an adjacency tensor 
$\mathcal{A}$ is defined to shed light on the communication in a multiplex. The 
sensitivity of the transmission of information to perturbations of the entries of 
$\mathcal{A}$ is investigated, and indicates both the edges of the multiplex that can be 
removed and the edges that should be strengthened. 

\section*{Acknowledgment}
The authors would like to thank a referee for comments that improved the presentation.
Research by SN was partially supported by a grant from SAPIENZA Universit\`a di Roma 
and by INdAM-GNCS.

\section*{Data availability statement} 
Data sharing not is applicable to this article as no new datasets were generated during the 
current study.

\section*{Conflict of interest} 
The authors declare that they have no conflict of interest. 

\section*{Funding}
See acknowledgement. 

\section*{Authors' contributions}
All authors contributed equally to the paper.

\section*{Ethical Approval}
Not Applicable.

\bibliographystyle{plain}

\end{document}